\documentclass{amsart}

\usepackage[utf8]{inputenc}
\usepackage{amsmath,amssymb,amsthm,units, stmaryrd,stackrel,relsize,bm}

\newtheorem{theorem}{Theorem}[section]

\newtheorem{lemma}[theorem]{Lemma}
\newtheorem{sublemma}[theorem]{Sublemma}

\newtheorem{corollary}[theorem]{Corollary}

\theoremstyle{definition}
\newtheorem{definition}[theorem]{Definition}

\theoremstyle{remark}
\newtheorem{remark}[theorem]{Remark}

\usepackage[
    colorlinks=true, 
    citecolor=red,
    linkcolor=blue,
    urlcolor=blue]{hyperref}
\newcommand\ad{{\mathsf{AD}}}

\newcommand\zfc{{\mathsf{ZFC}}}
\newcommand\zf{{\mathsf{ZF}}}
\newcommand\kp{{\mathsf{KP}}}

\newcommand\dc{{\mathsf{DC}}}

\newcommand\Ord{{\mathsf{Ord}}}

\newcommand\lth{\textsc{lth}}

\newcommand\ind{{\mathsf{IND}}}
\newcommand\coind{{\mathsf{coIND}}}

\newcommand\Def{\text{Def}}

\newcommand{\PI}{\boldsymbol\Pi}
\newcommand{\SIGMA}{\boldsymbol\Sigma}
\newcommand{\DELTA}{\boldsymbol\Delta}

\newcommand\quotesleft{``}

\makeatletter
\def\Ddots{\mathinner{\mkern1mu\raise\p@
\vbox{\kern7\p@\hbox{.}}\mkern2mu
\raise4\p@\hbox{.}\mkern2mu\raise7\p@\hbox{.}\mkern1mu}}
\makeatother

\begin{document}
\title{$F_\sigma$ Games and Reflection in $L(\mathbb{R})$}
\author{J. P. Aguilera}
\address{Institute of Discrete Mathematics and Geometry, Vienna University of Technology. Wiedner Hauptstra{\ss}e 8--10, 1040 Vienna, Austria.}
\email{aguilera@logic.at}

\begin{abstract}
We characterize the determinacy of $F_\sigma$ games of length $\omega^2$ in terms of determinacy assertions for short games. Specifically, we show that $F_\sigma$ games of length $\omega^2$ are determined if, and only if, there is a transitive model of $\kp+\ad$ containing $\mathbb{R}$ and reflecting $\Pi_1$ facts about the next admissible set.

As a consequence, one obtains that, over the base theory $\kp + \dc + \quotesleft \mathbb{R}$ exists,'' determinacy for $F_\sigma$ games of length $\omega^2$ is stronger than $\ad$, but weaker than $\ad + \Sigma_1$-separation.
\end{abstract}
\date{\today $\,$ (compiled)}
\clearpage
\maketitle

\setcounter{tocdepth}{1}
\tableofcontents
\section{Introduction}
We study the consistency strength of $F_\sigma$-determinacy for games of length $\omega^2$. We see that the situation here is tied with that of short games. For example, over the base theory 
\[\text{$\kp + \dc $ + \quotesleft $\mathbb{R}$ exists,''}\] 
$F_\sigma$-determinacy for games of length $\omega^2$ is much stronger than $\ad$; however, it is much weaker than $\ad$ + $\Sigma_1$-separation. Over $\zfc$, $F_\sigma$-determinacy for games of length $\omega^2$ is stronger than the existence of a transitive model of $\kp$ + $\ad$ containing all reals; yet weaker than the existence of a transitive model of $\kp$ + $\ad$ + $\Sigma_1$-separation containing all reals.
The consistency strength of this theory is hard to describe in terms of large cardinals: determinacy assumptions for games of length $\omega^2$ (over $\zfc$ or over $\kp + \dc$ + \quotesleft $\mathbb{R}$ exists'') for all pointclasses between the clopen sets and the Borel sets lie strictly between the existence of all finite amounts of Woodin cardinals (as a schema) and the existence of infinitely many Woodin cardinals, in terms of consistency strength (see \cite{AgThesis}).

 Much like in the case of short $F_\sigma$ games, long games are tied to reflection for set-theoretic formulae of a certain complexity, except that---rather than in $L$---one needs to consider reflection in the $L(\mathbb{R})$-hierarchy. Given an admissible set, let $A^+$ denote the \emph{next admissible set}, in the sense of Barwise-Gandy-Moschovakis \cite{BGM}; i.e., $A^+$ is the smallest admissible set containing $A$. By convention, all admissible sets are transitive.

\begin{definition}
An admissible set $A$ is \emph{$\Pi^+_1$-reflecting} if for every $\Pi_1$ formula $\phi$ with parameters in $A$ and one free variable, if $A^+\models\phi(A)$, then there is some $B \in A$ with all parameters in $\phi$ such that $B^+\models\phi(B)$. 
\end{definition}

Our main theorem is:
\begin{theorem}\label{mainintrofsigma}
The following are equivalent:
\begin{enumerate}
\item $\SIGMA^0_2$-determinacy for games of length $\omega^2$;
\item There is a $\Pi^+_1$-reflecting model of $\ad$ containing $\mathbb{R}$.
\end{enumerate}
\end{theorem}

Precursors to this work include, on the side of long games, Blass' \cite{Bl75} theorem that determinacy for all games of length $\omega^2$ is equivalent to $\ad_\mathbb{R}$, Neeman's \cite{Ne04} work on games of countable length, Trang's \cite{Tr13} work on analytic games of additively indecomposable length, as well as previous work on games of length $\omega^2$ that are open (see Theorem \ref{TheoremKPOpen} below), clopen \cite{Ag18b}, Borel \cite{Ag18c}, or projective \cite{AgMu}. On the side of $F_\sigma$ games we mention Solovay's and Tanaka's \cite{Ta91} work in the contexts of subsystems of set theory and analysis, respectively, and Wolfe's \cite{Wo55} original proof of $F_\sigma$-determinacy.

\section{Preliminaries}
In this section, we collect some preliminary definitions and results. 
\subsection{Games of transfinite length}
We study two-player, perfect information games of length $\omega^2$ which are $F_\sigma$, i.e., $\SIGMA^0_2$-definable. These are games in which, given $A\in \mathcal{P}(\mathbb{R})\cap \SIGMA^0_2$, two players, Player I and Player II, alternate $\omega^2$-many turns playing natural numbers, thus producing a sequence $x \in \mathbb{N}^{\omega^2}$. Since the spaces $\mathbb{R}\setminus \mathbb{Q}$ and $\mathbb{N}^{\omega^2}$ are recursively homeomorphic, the sequence $x$ may be identified with a(n irrational) real number. Player I wins if $x \in A$; otherwise, Player II wins.

As above, we will often identify sequences of natural numbers of length $\omega^2$ with $\omega$-sequences of reals and with single reals. We will denote by 
\[x = \langle x_0, x_1, \hdots \rangle\]
the \emph{single} real number coding the sequence $(x_0, x_1,\hdots)$, via some fixed recursive coding. The precise coding used will be immaterial, except for the continuity property that the first $n$ digits $x$ should only depend on $x_0, x_1,\hdots, x_n$. Similarly, if $s$ and $t$ are infinite sequences, or reals coding them, we will denote by 
\[s^\frown t\]
the single real number coding the result of concatenating $s$ and $t$.

\subsection{Reflection in $L(\mathbb{R})$}
We recall the definition of the $L(\mathbb{R})$ hierarchy: $L_0(\mathbb{R})$ is defined to be $V_{\omega+1}$, the collection of all sets all of whose elements are hereditarily finite. $L_{\alpha+1}(\mathbb{R})$ is the set of all subsets of $L_\alpha(\mathbb{R})$ definable over $L_\alpha(\mathbb{R})$ with parameters. At limit stages, $L_\lambda(\mathbb{R})$ is the union of all $L_\alpha(\mathbb{R})$, for $\alpha<\lambda$. The function
$$\alpha \mapsto L_\alpha(\mathbb{R})$$
is $\Delta_1$ with $\mathbb{R}$ as a parameter.
The first step towards proving Theorem \ref{mainintrofsigma} is the observation that it suffices to restrict to admissible sets of the form $L_\alpha(\mathbb{R})$. 

We will use the following non-standard notation:
\begin{definition} \label{DefinitionPlusLR}
Let $\alpha$ be an ordinal. We denote by 
$$\alpha^+$$
the least ordinal $\beta$ such that $\alpha<\beta$ and $L_{\beta}(\mathbb{R})$ is admissible.
\end{definition}

\begin{definition}
An ordinal $\alpha$ is \emph{$\Pi^+_1$-reflecting} if $L_\alpha(\mathbb{R})$ is $\Pi^+_1$-reflecting.
\end{definition}

\begin{lemma}\label{LemmaPi+reflectionOrdinals}
Suppose $A$ is a $\Pi^+_1$-reflecting set and $\mathbb{R} \in A$. Let $\alpha = \Ord\cap A$; then, $\alpha$ is $\Pi^+_1$-reflecting.
\end{lemma}
\proof
Let $\phi$ be a $\Pi_1$ formula with parameters in $L_\alpha(\mathbb{R})$ such that 
$$\left(L_{\alpha}(\mathbb{R})\right)^+\models\phi(L_\alpha(\mathbb{R})).$$
Clearly, $\left(L_{\alpha}(\mathbb{R})\right)^+ = L_{\alpha^+}(\mathbb{R})$. Moreover, $A^+$ is an admissible set containing $L_\alpha(\mathbb{R})$, so $L_{\alpha^+}(\mathbb{R})\subset A^+$.
Let $\psi$ be the formula in the language of set theory with parameters in $A \cup \{A\}$ asserting that for all $\beta$, if no $\gamma \in (\alpha,\beta)$ is admissible, then $L_\beta(\mathbb{R})\models \phi(L_\alpha(\mathbb{R}))$. Being admissible is expressible internally, $\alpha$ is $\Delta_1$-definable from $A$, and---as remarked earlier---$L_\beta(\mathbb{R})$ is $\Delta_1$-definable from $\beta$ and $\mathbb{R}$ (which belong to $L_\alpha(\mathbb{R})$). Hence, $\psi$ is $\Pi_1$. Moreover, $A^+ \models \psi(A)$, so, by reflection, there is $B \in A$ such that $B^+\models\psi(B)$ and $B$ contains all parameters in $\psi$. In particular, $\mathbb{R} \in B$.
Let $\beta = \Ord\cap B$; since $B \in A$, $\beta \in A$, so, in particular, $\beta<\alpha$. $L_{\beta^+}(\mathbb{R})\subset B^+$. By choice of $\psi$, we obtain that $L_{\beta^+}(\mathbb{R})\models \phi(L_\beta(\mathbb{R}))$, as desired.
\endproof

\begin{corollary}\label{CorollaryCharacPiplusRef}
The following are equivalent:
\begin{enumerate}
\item There is a $\Pi^+_1$-reflecting model of $\ad$ containing $\mathbb{R}$; and
\item There is a $\Pi^+_1$-reflecting ordinal $\alpha$ such that $L_\alpha(\mathbb{R})\models\ad$.
\end{enumerate}
\end{corollary}
\proof
Suppose that $A$ is a $\Pi^+_1$-reflecting model of $\ad$ containing $\mathbb{R}$ and let $\alpha = \Ord\cap A$. By Lemma \ref{LemmaPi+reflectionOrdinals}, $L_\alpha(\mathbb{R})$ is $\Pi^+_1$-reflecting. Moreover, it is a subset of $A$ containing all reals, and thus all possible strategies for games; hence, $L_\alpha(\mathbb{R})\models\ad$. The converse is immediate.
\endproof

This motivates the following definition:
\begin{definition}
We denote by $\sigma_\mathbb{R}$ the least $\Pi^+_1$-reflecting ordinal.
\end{definition}
$\sigma_\mathbb{R}$ is the analog of the ordinal $\sigma^1_1$ for $L(\mathbb{R})$. Unlike $\sigma^1_1$, the set-theoretic properties of $\sigma_\mathbb{R}$ are not decidable within $\zf$. For instance,
under $\ad$, $\sigma_\mathbb{R}$ is big---it is a limit of weakly inaccessible cardinals. In contrast, if $V = L$, then $\sigma_\mathbb{R}<\omega_2$.

\begin{definition}
Let $A\subset\mathbb{R}\times\mathbb{R}$ and write $A_x = \{y\in\mathbb{R}: (x,y)\in\mathbb{R}\}$. We write $\Game^\mathbb{R} A$ for the set of all $x\in\mathbb{R}$ such that Player I has a winning strategy for the game on $\mathbb{R}$ with payoff $A_x$. We define
$$\Game^\mathbb{R}\SIGMA^0_2 = \{\Game^\mathbb{R} A : A \in \SIGMA^0_2\}.$$
\end{definition}

(Classes such as $\Game^\mathbb{R}\SIGMA^0_1$ or $\Game^\mathbb{R}\PI^0_1$ are defined analogously.)
To prove Theorem \ref{mainintrofsigma}, it suffices to prove that determinacy for $F_\sigma$ games of length $\omega^2$ is equivalent to the fact that $L_{\sigma_\mathbb{R}}(\mathbb{R})\models\ad$. Following previous proofs of determinacy for games of length $\omega^2$, our first goal is to locate winning strategies for $\SIGMA^0_2$ games on $\mathbb{R}$. We shall show that
\[\mathcal{P}(\mathbb{R}) \cap \Sigma_1^{L_{\sigma_\mathbb{R}}(\mathbb{R})}= \Game^\mathbb{R}\SIGMA^0_2.\]

The main tool for this is the theory of inductive definitions on $\mathbb{R}$, the basics of which we recall next.

\subsection{Inductive definitions}
Suppose that 
$\phi:\mathcal{P}(\mathbb{R}^m)\to\mathcal{P}(\mathbb{R}^m)$
is an operator. We define sets $\phi^\lambda\subset\mathbb{R}^m$ inductively by 
\begin{align*}
\phi^0 &= \varnothing,\\
\phi^{<\lambda} &= \bigcup_{\mu<\lambda}\phi(\phi^{<\mu}),\\
\phi^{\lambda} &= \phi^{<\lambda}\cup \phi(\phi^{<\lambda}),\\
\phi^\infty &=\bigcup_{\lambda\in\Ord}\phi^{\lambda}.
\end{align*}

The least $\kappa$ such that $\phi^\infty = \phi^\kappa$ is called the \emph{closure ordinal} of $\phi$ and denoted $|\phi|$. If $\Gamma$ is a class of operators, we write 
\index{$|\Gamma|$}
$$|\Gamma| = \sup\{|\phi|:\phi \in \Gamma\}.$$

\begin{definition}\label{DefGamma-inductive}
Let $\Gamma$ be a class of operators
$$\phi:\mathcal{P}(\mathbb{R}^m)\to\mathcal{P}(\mathbb{R}^m).$$
We say that $R\subset\mathbb{R}^n$ is \emph{$\Gamma$-inductive} \index{inductive (relation)} if there is $b \in\mathbb{R}^{m-n}$ such that for all $a \in \mathbb{R}^n$,
$$a \in R \text{ if, and only if, } (a,b) \in \phi^\infty.$$
\end{definition}
An operator $\phi$ is defined by a formula $\psi(x,X)$ if, and only if, 
$$\phi(X) = \{a \in \mathbb{R}^m: \psi(a,X)\}.$$
Thus, it is natural to consider classes of operators specified in terms of definability. Let $\Gamma$ be a pointclass; we say that an operator $\phi$ is in $\Gamma$ if it is defined by a formula in $\Gamma$ with an additional predicate symbol $X$. We say that an operator is \emph{positive} \index{positive operator} if this additional predicate symbol appears only positively, i.e., not in the scope of any negations (or in the antecedent of implications).
An operator is \emph{monotone}\index{monotone operator} if $X\subset Y\subset\mathbb{R}$ implies $\phi(X)\subset\phi(Y)$. Every positive operator is monotone, of course. We may also speak of a formula defining an operator being \emph{positive} under the obvious circumstances.

\begin{definition}
Let $\Gamma$ be a pointclass.
\begin{enumerate}
\item $\Gamma{-}\ind$ is the pointclass of all $\Gamma$-inductive sets. 
\item $\Gamma^{pos}{-}\ind$ is the pointclass of all (positive-$\Gamma$)-inductive sets. 
\item $\Gamma^{mon}{-}\ind$ is the pointclass of all (monotone-$\Gamma$)-inductive sets.
\end{enumerate}
\end{definition}

\begin{definition}
A subset of $\mathbb{R}$ is \emph{inductive} if it is positive analytical inductive, i.e., letting $\ind$ be the class of all inductive sets, we have
$$\ind = \bigcup_{n\in\mathbb{N}}(\SIGMA^{1,pos}_n{-}\ind).$$
A subset of $\mathbb{R}$ is \emph{co-inductive} if its complement is inductive. We write $\coind$ for the class of co-inductive sets.
\index{inductive (pointclass)}
\end{definition}

Pointclasses of the form $\SIGMA^1_n$ (or of other forms) can be relativized by allowing sets of reals as parameters. This leads to a relativized version of the inductive sets:
\begin{definition}
Let $X\subset\mathbb{R}$. A subset of $\mathbb{R}$ is \emph{inductive on $X$} if it is positive analytical-on-$X$ inductive, i.e., letting $\ind(X)$ be the class of all sets inductive on $X$, we have
$$\ind(X) = \bigcup_{n\in\mathbb{N}}(\SIGMA^{1,pos}_n(X){-}\ind).$$
\end{definition}
In the definition above, ``positive'' refers to the variable on which the induction is carried out, not to the parameter $X$.

Let $\kappa$ be the closure ordinal of the inductive sets, i.e., the supremum of closure ordinals of positive analytical inductive definitions. Then $L_\kappa(\mathbb{R})$ is the smallest model of $\kp$ containing $\mathbb{R}$ and 
$$\ind = \mathcal{P}(\mathbb{R})\cap\Sigma_1^{L_\kappa(\mathbb{R})},$$
where $\Sigma_1$-definability in the equation above allows for parameters.
Subsets of $\mathbb{R}$ that are both inductive and co-inductive are called \emph{hyperprojective}. By $\Delta_1$-separation, these are the sets of reals in $L_\kappa(\mathbb{R})$. 

In general, let $\alpha$ be an ordinal and suppose that $\alpha$ is smaller than the least non-$\mathbb{R}$-projectible ordinal, i.e., that there is a surjection
$$\rho:\mathbb{R}\to L_\alpha(\mathbb{R})$$
which is $\Sigma_1$-definable over $L_\alpha(\mathbb{R})$ (this restriction can be relaxed, but all ordinals relevant to this article will be of this form). Let $X$ be a set of reals coding $L_\alpha(\mathbb{R})$ this way. Then, if $\kappa_X$ is the closure ordinal of the $X$-inductive sets, we have that $L_{\kappa_X}(\mathbb{R})$ is the smallest model of $\kp$ containing $L_\alpha(\mathbb{R})$ and 
$$\ind(X) = \mathcal{P}(\mathbb{R})\cap\Sigma_1^{L_{\kappa_X}(\mathbb{R})}.$$

An alternate definition of the inductive sets is given by the following theorem:
\begin{theorem}[Moschovakis \cite{Mo72}]
The following are equivalent:
\begin{enumerate}
\item $A\subset\mathbb{R}$ is inductive,
\item there is a projective (or even analytic) set $B\subset\mathbb{R}$ such that for all $x\in\mathbb{R}$,
\[x \in A \leftrightarrow \exists x_0\,\forall x_1\,\exists x_2\,\hdots \exists n\, (x,\langle x_0, \hdots, x_n\rangle) \in B.\]
\end{enumerate}
\end{theorem}

Moschovakis' theorem has the following easy consequence:
\begin{corollary}\label{TheoremCharacIndXGames}
$\ind(X) = \Game^\mathbb{R}\SIGMA^0_1(X)$.
\end{corollary}
\proof
Assume $X = \varnothing$ for simplicity. To see that 
\[\ind \subset \Game^\mathbb{R}\SIGMA^0_1,\]
it suffices to notice that a statement of the form
\[\exists x_0\,\forall x_1\,\exists x_2\,\hdots \exists n\, (x,\langle x_0, \hdots, x_n\rangle) \in B\]
can be decided by an open (in our sense) game on $\mathbb{R}$ (the details can be verified e.g., as in \cite{Ag18a}). Conversely, let $U \in \SIGMA^0_1\cap\mathcal{P}(\mathbb{R}^2)$. 
Write 
$$U = \bigcup_{i\in\mathbb{N}}U_i$$
as a countable union of basic clopen sets.
Thus, 
$$x \in \Game^\mathbb{R}U \text{ if, and only if, } \exists x_0\in\mathbb{R}\,\forall x_1\in\mathbb{R}\,\hdots\,\exists n\in\mathbb{N}\, (x,\langle x_0,x_1,\hdots\rangle) \in U_n.$$
Since $U_n$ is basic open, it is of the form
$$O(s,t) = \{(x,y)\in\mathbb{R}: (s,t)\sqsubset (x,y)\}.$$
Let $k_n$ be the length of $s$ and $t$. By our choice of coding of infinite sequences of real numbers by real numbers, the first $k_n$ digits of a sequence $\langle x_0,x_1,\hdots\rangle$ depend only on $x_0,\hdots, x_{k_n}$.
Clearly, the set
\begin{align*}
U_n^* &= \big\{(x, \langle k, x_0,\hdots, x_{k}\rangle) : \text{ $k$ is greater than or equal to the length of the }\\
&\quad \text{ unique finite  sequences $s$ and $t$ such that $U_n = O(s,t)$, and }\\
&\quad \quad \exists x_{k+1}\,\exists x_{k+2}\,\hdots\,(x,\langle x_0,x_1,\hdots\rangle) \in U_n)\big\};
\end{align*}
is projective. Let
\begin{align*}
A = \big\{x\in\mathbb{R}: \exists x_0\in\mathbb{R}\,\forall x_1\in\mathbb{R}\,\hdots\,\exists n, m\in\mathbb{N}\, (x,\langle n, x_0,x_1,\hdots, x_n\rangle) \in U^*_m\big\}.
\end{align*}
We claim that $A$ is inductive. To see this, first notice that the collection of all pairs $(x,\langle n, x_0,x_1,\hdots, x_n\rangle)$ belonging to some $U^*_m$ is projective. Moreover, the real $\langle n, x_0,x_1,\hdots, x_{\langle n, m\rangle}\rangle$ coding the tuple $(n, x_0,x_1,\hdots, x_{\langle n, m\rangle})$ agrees on the first $n$ digits with the real $\langle n, x_0,x_1,\hdots, x_{n}\rangle$ coding the tuple $(n, x_0,x_1,\hdots, x_n)$. Thus, it follows from the definition of $U_m^*$ that
$$(x,\langle n, x_0,x_1,\hdots, x_n\rangle) \in U^*_m$$
implies
$$(x,\langle \langle n, m\rangle, x_0,x_1,\hdots, x_{\langle n, m\rangle }\rangle) \in U^*_m.$$
Hence, $A$ is inductive. Finally, writing $k_n$ for the length of the finite sequences $s$ and $t$ such that $U_n = O(s,t)$,
\begin{align*}
x\in A\text{ iff }
&\exists x_0\in\mathbb{R}\,\forall x_1\in\mathbb{R}\,\hdots\,\exists n,m\in\mathbb{N}\, (x,\langle n, x_0,x_1,\hdots, x_n\rangle) \in U^*_m\\
\text{ iff }&\exists x_0\in\mathbb{R}\,\forall x_1\in\mathbb{R}\,\hdots\,\exists n,m\in\mathbb{N}\ k_m\leq n\ \wedge \\
&\qquad\exists y_{n+1}\,\exists y_{n+2}\,\hdots\,(x,\langle x_0,x_1,\hdots, x_n, y_{n+1}, \hdots\rangle) \in U_m\\
\text{ iff }&\exists x_0\in\mathbb{R}\,\forall x_1\in\mathbb{R}\,\hdots\,\exists m\in\mathbb{N}\, (x,\langle x_0,x_1,\hdots\rangle) \in U_m
\end{align*}
Here, the first and second equivalences hold by definition. The third equivalence holds because whether $(x,\langle x_0,x_1,\hdots\rangle)$ belongs to $U_m$ depends only on the first $k_m$ digits of $\langle x_0,x_1,\hdots\rangle$ (by definition of $k_m$), which in turn depend only on $x_0, x_1, \hdots, x_{k_m}$.
\endproof

We mention some results on the relation between inductive definitions, monotone inductive definitions, and positive inductive definitions.\footnote{We would like to thank A. S. Kechris for sharing the notes \cite{HK76h} with us.}
\begin{theorem}[Harrington-Kechris \cite{HK76,HK76h}] \label{HarringtonKechrisTheorem}
Let $\mathcal{F}$ be a collection of operators containing $\Pi^1_1$ and closed under $\wedge,\vee,\exists^\mathbb{R}$ and recursive substitutions. Suppose that
\begin{enumerate}
\item $\text{WF}\in \breve{\mathcal{F}}$, and
\item $\breve{\mathcal{F}} \in \mathcal{F}^{mon}{-}\ind$.
\end{enumerate}
Then $\mathcal{F}^{mon}{-}\ind = \mathcal{F}{-}\ind$.
\end{theorem}

We state the following theorem without defining all notions involved, and afterwards state the instance in which we will be interested:
\begin{theorem}[Harrington-Moschovakis \cite{HM75}]
Let $Q$ be a quantifier on $\mathbb{R}$. Let $\Gamma$ be the pointclass of all sets which are positive $Q$-inductive and let $\Gamma'$ be the pointclass dual to $\Gamma$. Then,
\[\Gamma'{-}\ind = \Gamma'^{pos}{-}\ind.\]
\end{theorem}
In particular, letting $Q$ be the quantifier $\exists^\mathbb{R}$, we have $\Gamma = \ind$ and $\Gamma' = \coind$, and thus:

\begin{theorem}[Harrington-Moschovakis]\label{Harrington-MoschovakisTheorem}
$$\coind{-}\ind = \coind^{pos}{-}\ind.$$
\end{theorem}

We mention that Harrington and Moschovakis' theorem as stated in \cite{HM75} is more general and e.g., implies the following classical result:
\begin{theorem}[Grilliot]\label{GrilliotsTheorem}
$\Sigma^1_1{-}\ind = \Sigma_1^{1,pos}{-}\ind.$
\end{theorem}

We finish by recalling the definition of a Spector class (on $\mathbb{R}$).

\begin{definition}
A pointclass $\Gamma$ is called\footnote{This differs from the notion of \quotesleft Spector \emph{pointclass}'' in Moschovakis \cite{Mo09}.} a \emph{Spector class on $\mathbb{R}$} if the following conditions hold:
\begin{enumerate}
\item $\Gamma$ is closed under $\wedge, \vee, \exists^\mathbb{R}, \forall^\mathbb{R}$;
\item $\Gamma$ contains all analytical relations;
\item $\Gamma$ is parametrized by $\mathbb{R}$;
\item $\Gamma$ has the prewellordering property.
\end{enumerate}
\index{Spector class (on $\mathbb{R}$)}
\end{definition}

\subsection{Short $F_\sigma$ games}
In order to locate winning strategies for $\SIGMA^0_2$ games on $\mathbb{R}$ within the $L(\mathbb{R})$-hierarchy, it might be helpful to recall the location of winning strategies for $\Sigma^0_2$ games on $\mathbb{N}$ in the $L$-hierarchy. This is a result of Solovay and is obtained by combining three theorems in recursion theory. 

The first one of them is also due to Solovay:
\begin{theorem}[Solovay]\label{SolovaysTheoremFsigma}
$\Sigma^{1,pos}_1{-}\ind = \Game\Sigma^0_2.$
\end{theorem}
A proof of Solovay's theorem can be found in \cite{Ke16} or in \cite{Mo09}. It will adapt to prove the analog for games on $\mathbb{R}$ later, as well as the fact that $\Game^\mathbb{R}\SIGMA^0_2$-determinacy implies determinacy for $\SIGMA^0_2$ games of length $\omega^2$.
The second theorem is Grilliot's Theorem \ref{GrilliotsTheorem} which, together with Theorem \ref{SolovaysTheoremFsigma} implies that
$$\Sigma^{1}_1{-}\ind = \Game\Sigma^0_2.$$

The final theorem is due to Aczel and Richter. Recall that $\sigma^1_1$ denotes the least $\Sigma^1_1$-reflecting  ordinal.
\begin{theorem}[Aczel-Richter \cite{AcRi74}]\label{AczelRichter}
$|\Sigma^1_1| = \sigma^1_1$, where $|\Sigma^1_1|$ refers to closure under inductive operators on $\mathbb{N}$.
\end{theorem}

Combining everything, one sees that, in order to know which player wins a (lightface) $\Sigma^0_2$ game on $\mathbb{N}$, one needs not search beyond $L_{\sigma^1_1}$ for a winning strategy.
Our plan for locating Player I's winning strategies for $\SIGMA^0_2$ games on $\mathbb{R}$ will be to follow the same steps as in the argument for games on $\mathbb{N}$.

\section{Coinductive operators}
\begin{theorem}\label{LemmacoINDGames}
$\Game^\mathbb{R}\SIGMA^0_2 = \coind^{pos}{-}\ind$
\end{theorem}
\proof
We begin by quoting Wolfe's proof of $\Sigma^0_2$-determinacy and Solovay's proof of Theorem \ref{SolovaysTheoremFsigma} to prove that
\begin{equation}\label{LemmacoINDGames1}
\Game^\mathbb{R}\SIGMA^0_2 \subset \coind^{pos}{-}\ind.
\end{equation}
Although the argument for this inclusion is very similar to the one for games on $\mathbb{N}$, a variation of it will be used below in the proof of Lemma \ref{LemmaShorteningSigma02} below, so we prefer to include all the details.

Let $A\subset\mathbb{R}^2$ be a $\Sigma^0_2$ set (where we assume no parameter is needed; the general result follows by relativization). Say $A$ is given by 
$$(x,y)\in A \leftrightarrow \exists n\,\forall m\, (n,x\upharpoonright m,y\upharpoonright m) \in P,$$
where $P$ is recursive. 
The proof proceeds by showing that winning a certain game $G^*(\langle\rangle)$ equivalent to $A_x$ is equivalent to the membership of $x$ in a set in $\coind^{pos}{-}\ind$. For $s$ a finite sequence of real numbers of even length, define $G^*(s)$ to be the following game:
\begin{enumerate}
\item Players I and II alternate turns playing real numbers $\alpha_n$.
\item After infinitely many rounds have taken place, Player I wins if, letting 
\begin{align*}
t(m) = s^\frown\langle\alpha_0,\alpha_1,\hdots,\alpha_m\rangle
\end{align*} 
(i.e., letting $t$ be the real coding the result of concatenating $s$ with the first $m$ moves of the play $(\alpha_0,\alpha_1,\hdots)$), 
we have
$$\exists n\,\forall m\, \forall m'<m\, (n,x\upharpoonright m', t(m)\upharpoonright m') \in P.$$
\end{enumerate}
Note that by our conventions on coding sequences of reals by single reals, we have
\begin{align*}
s^\frown\langle\alpha_0,\alpha_1,\hdots\rangle\upharpoonright m 
= t(m)\upharpoonright m.
\end{align*} 
This implies that Player I has a winning strategy for the game on $A$ if, and only if, she has one for the game on $G^*(\langle\rangle)$. The game $G^*(s)$ is like $G^*(\langle\rangle)$, except that we assume that $s$ has already been played. We shall show that Player I having a winning strategy is in $\coind^{pos}{-}\ind$.

Let $s$ be a finite sequence of real numbers of even length and $X$ be a set of reals. Consider the following game, $G(X,s)$:
\begin{enumerate}
\item Players I and II alternate turns playing real numbers $\alpha_n$.
\item After infinitely many rounds have taken place, Player I wins if, and only if, letting 
$$t(2m) = s^\frown\langle\alpha_0,\hdots,\alpha_{2m}\rangle,$$ 
one of the following holds for each $m\in\mathbb{N}$:
\begin{enumerate}
\item there is $n\leq \lth(s)$ such that
we have
$$\forall m' < m\, (n, x\upharpoonright m', t(2m)\upharpoonright m') \in P.$$
\item $t(2m) \in X$.
\end{enumerate}
\end{enumerate}
The formula
\begin{equation}\label{eqLemmacoIndGamesWSGXs}
\phi(s,X) = \text{\quotesleft $s$ has even length and Player I has a winning strategy in $G(X,s)$''}
\end{equation}
is clearly in $\Game^\mathbb{R}\Pi^0_1$ with an additional predicate for $X$ and $X$ appears positively in it. Since $\Pi^0_1$ games on $\mathbb{R}$ are determined by the Gale-Stewart theorem \cite{GS53}, the dual pointclass of $\Game^\mathbb{R}\Pi^0_1$ is $\Game^\mathbb{R}\Sigma^0_1$. By Corollary \ref{TheoremCharacIndXGames}, \eqref{eqLemmacoIndGamesWSGXs} is in $\coind^{pos}$. The claim is now that
\begin{equation*}
s \in \phi^\infty \leftrightarrow \text{\quotesleft Player I has a winning strategy in $G^*(s)$''}
\end{equation*}

As in Solovay's proof, it is first shown that if $s \in \phi^\xi$, then Player I has a winning strategy in $G^*(s)$, by induction on $\xi$. 
Suppose that this holds for $\phi^{<\xi}$ and that $s \in \phi^\xi$, so that Player I has a winning strategy in $G(\phi^{<\xi},s)$, say, $\sigma$. A winning strategy for $G^*(s)$ is obtained as follows:
Player I begins by playing $G^*(s)$ by $\sigma$ so long as the first winning condition of $G(s,X)$ is satisfied, i.e., so long as after round $m$, letting $t(m) = s^\frown\langle\alpha_0,\alpha_1,\hdots,\alpha_m\rangle$, we have
$$\exists n\leq \lth(s)\, \forall m' < m\, (n, x\upharpoonright m', t(2m)\upharpoonright m') \in P.$$
If after some round this condition is not satisfied, then (since $\sigma$ is a winning strategy), we must have 
$$\exists \zeta<\xi\, \Big(s^\frown\langle\alpha_0,\hdots,\alpha_{2m}\rangle \in \phi^{\zeta}\Big),$$
in which case the induction hypothesis yields a winning strategy for the game $G^*(s^\frown\langle\alpha_0,\hdots,\alpha_{2m}\rangle)$ which should now be followed.

Conversely, if Player I has a winning strategy $\sigma$ in $G^*(s)$ then $s$ must belong to $\phi^\infty$, for otherwise (by monotonicity of $\phi$ and determinacy of closed games) Player II has a winning strategy $\tau$ in $G(s,\phi^\infty)$. If so, we could face off the strategies $\sigma$ and $\tau$ against each other. Since $\tau$ is winning for Player II, after finitely many rounds, one will reach a partial play $t(2m_1) = s^\frown\langle\alpha_0,\hdots,\alpha_{2m_1}\rangle$ such that the following hold:
\begin{enumerate}
\item $\forall n\leq \lth(s)\, \exists m' < m_1\, (n, x\upharpoonright m', t(2m_1)\upharpoonright m') \not\in P$; and
\item $t(2m_1) \not \in \phi^\infty$.
\end{enumerate}
The second condition yields a new strategy $\tau_1$ for Player II that can now be played against $\sigma$ in $G(t(2m_1), \phi^\infty)$ until some stage $m_2$ at which the two conditions above hold again, etc. Continuing this process infinitely often yields a play
$$t = s^\frown\langle\alpha_0,\alpha_1,\hdots\rangle$$
such that for each $k$,
$$\forall n\leq \lth(t(2m_k))\, \exists m' < m_k\, (n, x\upharpoonright m', t(2m_{k+1})\upharpoonright m') \not\in P,$$
so that
$$\forall n\,\exists m'(n, x\upharpoonright m', t\upharpoonright m') \not\in P,$$
contradicting the fact that $\sigma$ was a winning strategy for Player I in $G(s)$. This completes the proof of \eqref{LemmacoINDGames1}. 

The second step consists in showing that
\begin{equation}\label{LemmacoINDGames2}
\coind^{pos}{-}\ind \subset \Game^\mathbb{R}\SIGMA^0_2.
\end{equation}
Let $\phi(s,X)$ be a positive $\coind$ operator.

By a theorem of Kechris and Moschovakis (see \cite[Theorem 2.18]{Ke16}) and since $\SIGMA^0_2$ is  parametrized by $\mathbb{R}$ and has the prewellordering property, $\Game^\mathbb{R}\SIGMA^0_2$ is a Spector class on $\mathbb{R}$. By the \quotesleft Main Lemma'' of Moschovakis \cite{Mo74} (see also \cite[Theorem 1.7]{Ke16}), in order to see that $\phi^\infty \in \Game^\mathbb{R}\SIGMA^0_2$, it suffices to show that $\Game^\mathbb{R}\SIGMA^0_2$ is closed under $\phi$, i.e., that for all $A \in \Game^\mathbb{R}\SIGMA^0_2$, the set
$$A_\phi = \Bigg\{(x,z) : \phi\Big(x, \big\{y:(y,z) \in A\big\}\Big)\Bigg\}$$
is in $\Game^\mathbb{R}\SIGMA^0_2$. 
Let us put
$$A_z = \Big\{y\in\mathbb{R}:(y,z)\in A\Big\}$$
and verify that
$$\Big\{(x,z) : \phi\big(x, A_z\big)\Big\} \in \Game^\mathbb{R}\SIGMA^0_2.$$
Since $\phi$ is in $\coind$, it follows from Moschovakis \cite{Mo72} that there is a $\PI^1_1$ formula $\psi$ such that 
$$\phi(x,X) \leftrightarrow \exists y_0\,\forall y_1\, \hdots\, \forall n\, \psi(x,\langle y_0,\hdots, y_n\rangle,X),$$
say,
$$\phi(x,X) \leftrightarrow \exists y_0\,\forall y_1\, \hdots\, \forall n\, \forall w\, \psi_0(w, x,\langle y_0,\hdots, y_n\rangle,X),$$
for some arithmetical $\psi_0$.

To verify that $A_\phi$ is in $\Game^\mathbb{R}\SIGMA^0_2$, fix $x$ and $z$; we play the natural game given by the equation above:
\begin{enumerate}
\item Players I and II begin by playing real numbers $y_0, y_1,\hdots$ until Player II decides to move on to the next stage after turn $n$. If this never happens, Player I wins.
\item Player II plays $w \in\mathbb{R}$. Let $\theta_0 = \psi_0(w, x,\langle y_0,\hdots, y_n\rangle,A_z)$ and assume without loss of generality that $\theta_0$ has been rewritten without implications and without negations whose scope is not only an atomic formula.
\item If $\theta_k$ has been defined and is a formula in the language of second-order arithmetic with an additional predicate $A_z$, we proceed by cases:
\begin{enumerate}
\item If the outermost logical connective of $\theta_k$ is a disjunction, then Player I chooses one of the disjuncts; we set $\theta_{k+1}$ equal to this choice.
\item If the outermost logical connective of $\theta_k$ is a conjunction, then Player II chooses one of the conjuncts; we set $\theta_{k+1}$ equal to this choice.
\item If $\theta_k$ is atomic then either it is an atomic formula not involving $X$, in which case the game ends and Player I wins if, and only if, $\theta_k$ holds; or it is of the form $a \in A_z$ (since $\phi$ is positive). In the latter case, the game continues. 
\item Let $B \in \SIGMA^0_2$ be such that $A_z = \Game^\mathbb{R}B$. Players I and II alternate infinitely many turns playing reals $w_0, w_1, w_2,\hdots.$ At the end, Player I wins if, and only if, $(\langle w_0, w_1, \hdots\rangle, w, a, z) \in B$.
\end{enumerate}
\end{enumerate}
Clearly the winning condition for the game is $\SIGMA^0_2$. It is easy to verify that Player I has a winning strategy in this game if, and only if $(x,z) \in A_\phi$. This proves \eqref{LemmacoINDGames2}.
\endproof 

Together with Harrington and Moschovakis's Theorem \ref{Harrington-MoschovakisTheorem}, the previous result yields:

\begin{corollary}\label{corollaryCoInd=GSigma02}
$\coind{-}\ind = \Game^\mathbb{R}\SIGMA^0_2.$
\end{corollary}

\section{Coinductive operators and reflection}
Let $\preceq$ be a binary relation. As a convention, we will use $\prec$ to denote the strict part of $\preceq$ (i.e., $x\prec y$ whenever $x\preceq y$ and $y\not\preceq x$) and $\equiv$ to denote $\preceq$-equivalence (i.e., $x\equiv y$ whenever $x \preceq y$ and $y\preceq x$).  Similar conventions shall apply to variations of $\preceq$, e.g., by subscripts. Notice that any two of $\preceq$, $\prec$, and $\equiv$ determine the other one.

Given a binary relation $\preceq$ and an equivalence relation $E$ with the same field, $\preceq$ is said to commute with $E$ if $x\preceq y$, $xEa$, and $yEb$ imply $a\preceq b$, and similarly for the strict part $\prec$. A binary relation $\preceq$ \emph{commutes with itself} if it commutes with $\equiv$; this follows from transitivity. An induction on the complexity of formulae shows that if $\preceq$ is a binary relation that commutes with itself, then any first-order formula in the vocabulary $\{\prec, =\}$ holds of some tuple $(x_1,\hdots, x_n)$ of objects in $\text{field}(\preceq)$ if, and only if, it holds of any tuple $(y_1,\hdots, y_n)$ such that $x_i$ is $\equiv$-equivalent to $y_i$ for each $i$, so long as equality is interpreted as $\equiv$.

We define in the natural way the notion of a  preorder on $\mathbb{R}$ coding an initial segment of $L(\mathbb{R})$:
\begin{definition}
Let $\alpha$ be an ordinal. A \emph{coding} of $L_\alpha(\mathbb{R})$ is a preorder $\preceq$ on a subset of $\mathbb{R}$ such that
$$(L_\alpha(\mathbb{R}), \in, =) \cong
(\text{field}(\preceq), \prec, \equiv)/\equiv.$$
\end{definition}
The right-hand side of the displayed equation is the quotient structure of the structure  $(\text{field}(\preceq), \prec, \equiv)$ by $\equiv$. The relation $\equiv$ becomes equality in this quotient. If $\preceq$ is a coding of some $L_\alpha(\mathbb{R})$, then $\prec$ induces a wellfounded transitive relation on $\equiv$-equivalence classes and there is a surjection $\rho$ from field$(\preceq)$ to $L_{\alpha}(\mathbb{R})$ such that $x \prec y$ if, and only if, $\rho(x) \in \rho(y)$. If so then \emph{in this section only}, let us denote $\rho(x)$ by $|x|_\preceq$.

\begin{theorem}\label{theoremGeneralizedRichterAczel}
$|\coind| = \sigma_\mathbb{R}$.
\end{theorem}

The proof of Theorem \ref{theoremGeneralizedRichterAczel} will be divided into two lemmata:

\begin{lemma}\label{LemmaSigmaRReflects}
$|\coind|$ is $\Pi^+_1$-reflecting.
\end{lemma} 
\proof
The lemma will be proved by an argument similar to the one of \cite[Theorem 10.7]{AcRi74}. We have not checked whether the results of \cite{AcRi74} (e.g., the lemma proved in the appendix) hold in full in the context of recursion on $\mathbb{R}$ but, since we are only interested in proving Lemma \ref{LemmaSigmaRReflects} (and are not analizing, e.g., ordinals of the form $|\Sigma^1_n|$), we do not need to be very careful with issues of definability, which greatly simplifies the situation.

Given an inductive definition $\phi$, there is a natural way of associating to it a prewellordering on $\mathbb{R}$, namely, letting
\begin{align*}
w(x) =
\begin{cases}
\text{least $\xi$ such that $x\in\phi^{\xi+1}$}, &\text{ if it exists},\\
\infty, &\text{ otherwise};
\end{cases}
\end{align*}
we put
$$x\preceq_\phi y \text{ if, and only if, } w(x)\leq x(y).$$

Let $\phi$ be a universal coinductive operator, i.e., a coinductive operator such that whenever $\psi$ is a coinductive operator, we have
$$\psi(X) = \{x\in\mathbb{R}: (x,a) \in \phi(X)\}$$
for some $a\in\mathbb{R}$. (See Moschovakis \cite{Mo74} for a proof of the existence of universal coinductive sets; they can be defined uniformly in $X$.)
Below, we will define a coinductive operator $\Theta$ which simultaneously applies $\phi$, defines the prewellordering associated to $\phi$, and codes initial segments of $L_\alpha(\mathbb{R})$ along the prewellordering given by $\phi$. Given $X\subset\mathbb{R}$, write
\[(i,X) = \{(i,x):x\in\mathbb{R}\}.\]
Conversely, we write
$$X_i = \{x\in\mathbb{R}:(i,x) \in X\}.$$

Let $X \subset\mathbb{R}$ and $\preceq$ be a binary relation on $X$. Suppose moreover that $U$ is a set of real numbers and that $X$ consists only of tuples of real numbers none of whose first coordinate belongs to $U$. We define a set $\Def(X,\preceq)$ and a binary relation $\preceq^+$ such that if $\preceq$ has field $X$ and codes $L_{\alpha}(\mathbb{R})$, then $\preceq^+$ is a relation on $X \oplus\Def(X,\preceq)$  which codes $L_{\alpha+1}(\mathbb{R})$. The reason for mentioning $U$ at all will become clear soon; roughly, it contains indicators which we will use to identify which reals belong to $\Def(X,\preceq)$ and which do not.

$\Def(X,\preceq)$ is the set of all tuples $(x_0, \varphi, \vec b)$, where $x_0 \in U$, $\varphi$ is a formula of arity $\lth(\vec b) + 1$ in the vocabulary $\{=, \prec\}$ and $\vec b$ is a finite tuple of elements of $X$. For $x,y \in X\oplus \Def(X,\preceq)$, we put $x \prec^* y$ if, and only if, one of the following holds:
\begin{enumerate}
\item $x, y \in X$ and $x\prec y$;
\item $y \in \Def(X,\preceq)$ is of the form $(x_0, \varphi, \vec b)$ and $(X, \prec, \equiv) \models \varphi(\vec b, x)$.
\end{enumerate}
Afterwards, let $xEy$ if, and only if $x$ and $y$ have exactly the same $\prec^*$-predecessors, and put
$x \preceq^+ y$ if, and only if, $xEy$ or there are $aEx$ and $bEy$ such that $a \prec^* b$. Hence, the relation $\prec^*$ is a first approximation to $\prec^+$. It extends $\prec$ by specifying which of the elements of $X$ are smaller than which of the new elements. It may thus be that an element $x$ of $\Def(X,\preceq)$ turns out to have exactly the same $\prec^*$-predecessors as some other element $y$ of $\Def(X,\preceq)$ or of $X$. We would then like to force these $x$ and $y$ to be $\equiv^+$-equivalent, which is what we do. Note that this definition depended on $U$, so when in need of precision, we may write $\Def(X, \preceq, U)$ for $\Def(X, \preceq)$ and $\preceq^{+U}$ for $\preceq^+$. 
By inspecting the construction, one sees that both $\Def(X, \preceq, U)$ and $\preceq^+$ are, say, hyperprojective (in fact, much simpler) in $X$, $U$, and $\preceq$. Another feature of this definition is that if $\{\preceq_\iota\}_\iota$ is an increasing family of relations obtained this way, and $\preceq_0$ is a coding of some $L_{\alpha}(\mathbb{R})$ then $\bigcup_{\iota}\preceq_\iota$ is a coding of some $L_{\beta}(\mathbb{R})$ (this is  where the set $U$ in the definition is used).

We now define the operator $\Theta$. Given $X\subset\mathbb{R}$, $\Theta$ will map $X$ to a set $Y$ with four parts, $Y_0, Y_1, Y_2$, and $Y_3$ defined in terms of $X_0$, $X_1$, $X_2$, and $X_3$. If $X$ is of the right form, $Y_1$ will be a prewellordering of some length $\alpha+1$ and $Y_3$ will be a coding of $L_{\alpha+1}(\mathbb{R})$ with field $Y_2$.
\begin{enumerate}
\item $Y_0 = \phi(X_0)$;
\item $Y_1 = X_1 \cup \Big\{(x,y): x \in Y_0 \cup \text{field}(X_1) \text{ and } y\in Y_0\setminus \text{field}(X_1)\Big\}$;
\item $Y_2 = X_2 \cup \Def(X_2, X_3, U)$, where $U = Y_0\setminus \text{field}(X_1)$;
\item $Y_3 = X_3^{+U}$, where $U = Y_0\setminus \text{field}(X_1)$.
\end{enumerate}
The class of coinductive operators is closed under Boolean connectives and real quantifiers, so $\Theta$ is coinductive. The operator generates the inductive definition $\phi^\infty$ on its first component, while simultaneously coding $\preceq_\phi$ in its second component. At each stage $\xi$, all new elements added to the field of $\preceq_\phi$, i.e., those of $\preceq_\phi$-rank $\xi+1$, are used to define the new elements of the set $(\Theta^{\xi})_2$ which acts as the field of a relation $(\Theta^{\xi})_3$ that codes $L_{\xi+1}(\mathbb{R})$. 

Any fixed point of $\Theta$ will contain a fixed point of $\phi$ in its first component. Because $\phi$ is universal coinductive,
an argument as in \cite[Theorem 8.5]{AcRi74} (but using continuous reducibility in place of many-one reducibility\footnote{A proof of the Kleene's Recursion Theorem in this context can be found e.g., in \cite[Theorem 3.1]{KW10}.}) shows that 
$$|\Theta| = |\phi| = |\coind|.$$
Hence, we need to show that $|\Theta|$ is $\Pi^+_1$-reflecting (although, to show that there is a $\Pi^+_1$-reflecting ordinal ${\leq}|\coind|$---which is really the point of this lemma---, we do not need the argument from \cite[Theorem 8.5]{AcRi74}). Recall that, given a coding $\preceq$ of some $L_\alpha(\mathbb{R})$, in this section we write $|x|_\preceq$ for the element of $L_\alpha(\mathbb{R})$ coded by $x$, if any.

\begin{sublemma}
Let $\psi$ be a $\Pi_1$ formula in the language of set theory. Then, there is a coinductive operator $\Psi$, such that whenever $\lambda\leq |\Theta|$,
\begin{enumerate}
\item for all $c_1,\hdots, c_l \in (\Theta^\lambda)_2$,
\[L_{\lambda^+}(\mathbb{R})\models \psi\big(|c_1|_{(\Theta^\lambda)_3}, \hdots, |c_l|_{(\Theta^\lambda)_3}\big)\leftrightarrow (c_1,\hdots, c_l) \in \Psi(\Theta^\lambda).\]
\item for all $c_1,\hdots, c_l \in \mathbb{R}$, if there is some $i < l+1$ such that $c_i \not \in (\Theta^\lambda)_2$, then
\[(c_1,\hdots, c_l) \not\in \Psi(\Theta^\lambda).\]
\end{enumerate}
\end{sublemma}
\proof
The main observation is that by Barwise-Gandy-Moschovakis \cite[Lemma 2.9]{BGM}, for every admissible set $A$, a relation $P$ on $A$ is coinductive on $A$ with parameters if, and only if, it is $\Pi_1$ over $A^+$, with parameters in $A \cup \{A\}$. Moreover,   this correspondence is uniform, in that the definition of the coinductive relation depends only on the $\Pi_1$ formula and its parameters, and not on the admissible set $A$ (as long as it contains the parameters),\footnote{This can be shown directly for admissible sets of the form $L_\alpha(\mathbb{R})$: given such a set, one defines an operator, elementary in $L_\alpha(\mathbb{R})$, that successively outputs extensions of $L_{\alpha}(\mathbb{R})$, like $\Theta$ does. To ensure that the operator is elementary, one can have it e.g., add only $\Sigma_1$-definable sets at each stage, instead of all first-order--definable sets. One can ask the operator to also check at each limit stage whether the $\Pi_1$ fact with parameters in $L_{\alpha}(\mathbb{R})$ holds of the structure it outputs. Like we did with $\Theta$, this operator can be ensured to have closure ordinal $\alpha^+$ by asking it to generate a universal $\Sigma_1(L_\alpha(\mathbb{R}))$ inductive definition in parallel. The operator for that can be obtained e.g., from the universal $\Sigma_1(P)$ subset of $\mathbb{R}$, where $P$ is a subset of $\mathbb{R}$ coding $L_\alpha(\mathbb{R})$. The set $P$ can be obtained uniformly in a $\Sigma_1(L_\alpha(\mathbb{R}))$ way (for multiplicatively indecomposable $\alpha$; see Steel \cite[Lemma 1.4]{St08}).} and vice-versa.

Hence, given a $\Pi_1$-formula $\psi$, the set of all $a_1,\hdots,a_l \in L_{\lambda}(\mathbb{R})$ such that $$L_{\lambda^+}(\mathbb{R})\models \psi\big(a_1,\hdots,a_l \big)$$
is a coinductive subset of $L_{\lambda}(\mathbb{R})$ (and thus of $L_{\lambda+1}(\mathbb{R})$).
Since $(\Theta^\lambda)_3$ is a coding of $L_{\lambda+1}(\mathbb{R})$, there is a surjection
$$\rho:(\Theta^\lambda)_2 \to L_{\lambda+1}(\mathbb{R}),$$
such that a pair $(x,y)$ belongs to the strict part of $(\Theta^\lambda)_3$ if, and only if, $\rho(x) \in \rho(y)$.
It follows that the preimage of a coinductive subset of $L_{\lambda+1}(\mathbb{R})$ under $\rho$ is coinductive on $\mathbb{R}$ from the parameters $(\Theta^\lambda)_2$ and $(\Theta^\lambda)_3$. Hence, there is a coinductive operator $\Psi'$ such that 
\begin{enumerate}
\item for all $c_1,\hdots, c_l \in (\Theta^\lambda)_2$,
\[L_{\lambda^+}(\mathbb{R})\models \psi\big(|c_1|_{(\Theta^\lambda)_3}, \hdots, |c_l|_{(\Theta^\lambda)_3}\big)\leftrightarrow (c_1,\hdots, c_l) \in \Psi'\big((\Theta^\lambda)_2\oplus (\Theta^\lambda)_3\big).\]
\item for all $c_1,\hdots, c_l \in \mathbb{R}$, if there is some $i < l+1$ such that $c_i \not \in (\Theta^\lambda)_2$, then
\[(c_1,\hdots, c_l) \not\in \Psi'\big((\Theta^\lambda)_2\oplus (\Theta^\lambda)_3\big).\]
\end{enumerate}

Now, both $(\Theta^\lambda)_2$ and $(\Theta^\lambda)_3$ reduce to $\Theta^\lambda$ continuously, and this reduction is uniform in $\lambda$, so 
there is an operator $\Psi$ as desired.
\endproof

Let $\psi$ be a $\Pi_1$ sentence and let $\Psi$ be given by the sublemma, so that for all $\lambda\leq|\Theta|$ and all $c_1,\hdots, c_l \in (\Theta^\lambda)_2$, we have
\[L_{\lambda^+}(\mathbb{R})\models \psi\big(|c_1|_{(\Theta^\lambda)_3}, \hdots, |c_l|_{(\Theta^\lambda)_3}\big)\leftrightarrow (c_1,\hdots, c_l) \in \Psi(\Theta^\lambda).\]
Since $\phi$ was chosen to be a universal coinductive operator, $\Theta(X)$ is complete coinductive for all $X$, so there is a continuous function
$$g:\mathbb{R}^l\mapsto\mathbb{R}$$
such that for all $c_1,\hdots, c_l \in \mathbb{R}$ and all $X\subset\mathbb{R}$ do we have
\[(c_1,\hdots, c_l) \in \Psi(X) \leftrightarrow g(c_1,\hdots, c_l)\in \Theta(X). \]
Hence,
\begin{equation}\label{eqLemmaSigmaRReflects}
L_{\lambda^+}(\mathbb{R})\models \psi\big(|c_1|_{(\Theta^\lambda)_3}, \hdots, |c_l|_{(\Theta^\lambda)_3}\big)\leftrightarrow g(c_1,\hdots, c_l) \in \Theta(\Theta^\lambda).
\end{equation}
Suppose then that there are elements $\gamma_1,\hdots,\gamma_l$ of $L_{\lambda}(\mathbb{R})$ such that
$$L_{|\Theta|^+}(\mathbb{R})\models \psi\big(\gamma_1,\hdots, \gamma_l\big).$$
By construction, $(\Theta^{<\lambda})_3$ is a coding of $L_{\lambda}(\mathbb{R})$ for each $\lambda$, so $(\Theta^\infty)_3$ is a coding of $L_{|\Theta|}(\mathbb{R})$
with field $(\Theta^{\infty})_2$. Thus, there is a surjection 
$$\rho:(\Theta^{\infty})_2 \to L_{|\Theta|}(\mathbb{R})$$ 
such that the pair $(x,y)$ belongs to the strict part of $(\Theta^{\infty})_3$ if, and only if, $\rho(x) \in \rho(y)$. Pick $c_1,\hdots, c_l$ such that for each $0<i<l+1$, $\rho(c_i) = \gamma_i$. Then,
$$L_{|\Theta|^+}(\mathbb{R})\models \psi\big(|c_1|_{(\Theta^\infty)_3}, \hdots, |c_l|_{(\Theta^\infty)_3}\big).$$
By \eqref{eqLemmaSigmaRReflects},
$$g(c_1,\hdots, c_l) \in   \Theta(\Theta^\infty) \subset \Theta^\infty.$$ 
$|\Theta|$ is a limit ordinal, so there is some $\lambda<|\Theta|$ such that 
$$g(c_1,\hdots, c_l) \in   \Theta(\Theta^\lambda),$$
so that, by \eqref{eqLemmaSigmaRReflects},
$$L_{\lambda^+}(\mathbb{R})\models \psi\big(|c_1|_{(\Theta^\lambda)_3}, \hdots, |c_l|_{(\Theta^\lambda)_3} \big).$$ 
Notice that the sequence of sets $\{(\Theta^\lambda)_2:\lambda<|\Theta|\}$, as well as the sequence of isomorphisms witnessing that each $(\Theta^\lambda)_3$ is a coding of $L_{\lambda+1}(\mathbb{R})$ are strictly increasing, and this implies that 
$$|c_i|_{(\Theta^\lambda)_3} = |c_i|_{(\Theta^{\infty})_3}$$
for each $0<i<l+1$, which yields the result. This proves the lemma.
\endproof

\begin{lemma}
If $\kappa$ is $\Pi^+_1$-reflecting, then $|\coind|\leq\kappa$.
\end{lemma}
\proof
This proof is like that of \cite[Lemma 10.1]{AcRi74}.
Suppose $\kappa$ is $\Pi^+_1$-reflecting and let $\phi$ be a coinductive operator. Since $\kappa$ is $\Pi^+_1$-reflecting, it is  $\mathbb{R}$-recursively inaccessible. Because coinductive relations are universal relations on the next admissible set, it follows that if $A \in \mathcal{P}(\mathbb{R}) \cap L_\alpha(\mathbb{R})$, then $L_{\alpha^++1}(\mathbb{R})$ can compute coinductive relations on $A$.
Now, $L_\kappa(\mathbb{R})$ can compute the sequence $\{\phi^\lambda:\lambda<\kappa\}$ in a $\Pi_1$ way (with parameters): $x$ 
is the $\lambda$th element of this sequence if every transitive set of the form $L_\alpha(\mathbb{R})$ containing an increasing $(\lambda + 1)$-sequence of $\mathbb{R}$-admissible sets believes that $x$ is the $\lambda$th element of the sequence. Moreover, this definition is uniform on  $\mathbb{R}$-recursively inaccessible ordinals.
Thus, letting $\psi$ be a $\Pi_1$ formula such that
$$a \in \phi^{<\kappa} \leftrightarrow L_\kappa(\mathbb{R})\models\exists \lambda\,\psi(a,\lambda),$$
the question of whether $a$ belongs to $\phi(\phi^{<\kappa})$ is coinductive over $\phi^{<\kappa}$, hence over $L_\kappa(\mathbb{R})$, hence universal over $L_{\kappa^+}(\mathbb{R})$, so, by $\Pi^+_1$-reflection, $\phi(\phi^{<\kappa})\subset\phi^{<\kappa}$, as desired.
\endproof

This concludes the proof of Theorem \ref{theoremGeneralizedRichterAczel}. The main consequence of interest to us is:
\begin{corollary}\label{corollaryCharacterizationSigma02}
$\mathcal{P}(\mathbb{R}) \cap \Sigma_1^{L_{\sigma_\mathbb{R}}(\mathbb{R})} = \Game^\mathbb{R}\SIGMA^0_2$.
\end{corollary}
\proof
By Theorems \ref{LemmacoINDGames} and \ref{theoremGeneralizedRichterAczel},
we need to show that 
$$\mathcal{P}(\mathbb{R}) \cap \Sigma_1^{L_{|\coind|}(\mathbb{R})} = \coind{-}\ind.$$
That the class on the right-hand side is contained in the one on the left is clear. The converse follows from the proof of Lemma \ref{LemmaSigmaRReflects}: one can use a coinductive operator to inductively generate codes for initial segments of $L(\mathbb{R})$ below $|\coind|$, so an existential statement about $L_{|\coind|}(\mathbb{R})$ can be rephrased in terms of membership in the least fixed point of a coinductive operator.
\endproof

\section{$F_\sigma$-determinacy}
\begin{lemma}\label{LemmaShorteningSigma02}
Suppose that all sets in $\Game^\mathbb{R}\SIGMA^0_2$ are determined. Then, all $\SIGMA^0_2$ games of length $\omega^2$ are determined.
\end{lemma}
\proof
Note that $\SIGMA^0_2$ games of length $\omega$ with moves in $\mathbb{R}$ are determined.
Suppose that all sets in $\Game^\mathbb{R}\SIGMA^0_2$ are determined and let $A\in \SIGMA^0_2$. We show that the game of length $\omega^2$ with payoff $A$ is determined. Let us refer to this game as $G$. The proof begins very much like that of Theorem \ref{LemmacoINDGames}: we begin by replacing $G$ with a game $G^*(\langle\rangle)$. Without loss of generality we shall assume that there is a recursive set $P$ such that $A$ is given by 
$$x\in A \leftrightarrow \exists n\,\forall m\, (n,x\upharpoonright m) \in P,$$

For $s$ a finite sequence of real numbers of even length, define $G^*(s)$ to be the following game:
\begin{enumerate}
\item Players I and II alternate $\omega^2$-many turns playing natural numbers to produce a countable sequence $\{\alpha_n:n\in\mathbb{N}\}$ of real numbers.
\item After $\omega^2$-many rounds have taken place, Player I wins if, letting 
\begin{align*}
t(m) = s^\frown\langle\alpha_0,\alpha_1,\hdots,\alpha_m\rangle,
\end{align*} 
we have
$$\exists n\,\forall m\, \forall m'<m\, (n, t(m)\upharpoonright m') \in P.$$
\end{enumerate}
As before, Player I has a winning strategy for $G$ if, and only if, she has one for $G^*(\langle\rangle)$. Moreover, Player II has a winning strategy for $G$ if, and only if, she has one for $G^*(\langle\rangle)$. The game $G^*(s)$ is like $G^*(\langle\rangle)$, except that we assume that $s$ has already been played.

Let $n_s$ be a natural number, $s$ be a sequence of natural numbers of even length $\omega\cdot n_s$ and $X$ be a set of sequences of natural numbers of length ${<}\omega^2$. Consider the following game, $G(X,s)$:
\begin{enumerate}
\item Players I and II alternate $\omega^2$-many turns playing natural numbers to produce digits corresponding to a countable sequence $\{\alpha_n:n\in\mathbb{N}\}$ of real numbers.
\item After $\omega^2$-many rounds have taken place, Player I wins if, and only if, letting 
$$t(m) = s^\frown\langle\alpha_0,\hdots,\alpha_{m}\rangle,$$ 
one of the following holds for each $m\in\mathbb{N}$:
\begin{enumerate}
\item there is $n\leq n_s$ such that
we have
$$\forall m' < m\, (n, t(m)\upharpoonright m') \in P.$$
\item $t(m) \in X$.
\end{enumerate}
\end{enumerate}
Thus $G(X,s)$ is a game on $\mathbb{N}$ of length $\omega^2$. Write
\begin{equation}\label{eqLemmacoIndGamesWSGXsbis}
\phi(s,X) = \text{\quotesleft Player I has a winning strategy in $G(X,s)$.''}
\end{equation}
The argument from Theorem \ref{LemmacoINDGames} shows that if $s \in \phi^\infty$, then Player I has a winning strategy in $G^*(s)$. To complete the proof, it remains to show that if $s\not\in \phi^\infty$, then Player II has a winning strategy in $G^*(s)$.

Notice that $G(X,s)$ is a game of length $\omega^2$ with moves in $\mathbb{N}$ and payoff in $\Pi^0_1$ (using $X$ as a second-order oracle). We consider the following auxiliary game, which we shall denote by $H(X,s)$:
\begin{align*}
\begin{array}{c|ccccc}
I & &\sigma_0 & & \sigma_1 &\hdots\\
II &  \tau_0 & & \tau_1 & & \hdots
\end{array}
\end{align*}
Here, $\sigma_i$ and $\tau_i$ are real numbers coding strategies for games on $\mathbb{N}$ of length $\omega$. Player I wins the game if, and only if, the sequence
$$(s, \sigma_0*\tau_0*, \sigma_1*\tau_1,\hdots) $$
satisfies the winning condition of $G(X,s)$. Here, $\sigma*\tau$ denotes the result of facing off the strategies $\sigma$ and $\tau$ against each other. Thus, $H(X,s)$ is a game of length $\omega$ with moves in $\mathbb{R}$ and payoff in $\Pi^0_1(X)$.
We will need a sublemma:
\begin{sublemma}\label{SublemmaClopen}
Suppose that all subsets of $\mathbb{R}$ in $\Game^\mathbb{R}\SIGMA^0_1(X)$ are determined. Then, the following are equivalent:
\begin{enumerate}
\item Player I has a winning strategy in $H(X,s)$; and
\item Player I has a winning strategy in $G(X,s)$.
\end{enumerate}
Moreover, the games are determined.
\end{sublemma}
\proof[Proof Sketch]
The equivalence is proved by arguing as in the proof of \cite[Theorem 1.1]{Ag18b}. This is possible because the game $H(s,X)$ is $\Pi^0_1(X)$. The proof of \cite[Theorem 1.1]{Ag18b} uses determinacy of sets in $\Game^\mathbb{R}\SIGMA^0_1$; in this case, we need determinacy of sets in $\Game^\mathbb{R}\SIGMA^0_1(X)$, which we are assuming.
$H(X,s)$ is closed and thus determined by the Gale-Stewart theorem. The proof of \cite[Theorem 1.1]{Ag18b} also shows that if Player II has a winning strategy in $H(X,s)$, then she has one in $G(X,s)$, so that this game is also determined.
\endproof

\begin{remark}
Sublemma \ref{SublemmaClopen} essentially states that Lemma \ref{LemmaShorteningSigma02} holds if one replaces $\SIGMA^0_2$ by $\PI^0_1$.
The part of the proof of \cite[Theorem 1.1]{Ag18b} used in the proof of Sublemma \ref{SublemmaClopen} goes through for $\PI^0_1$ games. It would also go through for $\SIGMA^0_1$ games if one switched the order in which the players move in the definition of $H(X,s)$. The proof, however, does not go through for more complicated games. 
By different arguments (involving the theory of scales in $L(\mathbb{R})$), one can prove analogues of Lemma \ref{LemmaShorteningSigma02} for more complicated pointclasses, up to $\DELTA^1_1$ (see \cite{Ag18c}). The result of replacing $\SIGMA^0_2$ by $\PI^1_1$ in the statement of Lemma \ref{LemmaShorteningSigma02} is not provable in $\zfc$.
\end{remark}

It follows from Corollary \ref{TheoremCharacIndXGames} that if $X \in L_\alpha(\mathbb{R})$, then all winning strategies for games on $\mathbb{R}$ with payoff in $\SIGMA^0_1(X)$ are definable over $L_{\alpha^+}(\mathbb{R})$, so $\Game^\mathbb{R}\SIGMA^0_1(X)$ belongs e.g., to $L_{\alpha^{++}}(\mathbb{R})$.
\begin{sublemma}\label{SublemmaClosureInductive}
$\Game^\mathbb{R}\SIGMA^0_2$ is closed under the open-game-on-$\mathbb{R}$ quantifier; i.e.,
if $X \in \Game^\mathbb{R}\SIGMA^0_2$, then $\Game^\mathbb{R}\SIGMA^0_1(X)\subset\Game^\mathbb{R}\SIGMA^0_2$.
\end{sublemma}
\proof
This can be proved directly by an argument like the one of \eqref{LemmacoINDGames2} in Theorem \ref{LemmacoINDGames}. Alternatively, one can also appeal to
the \quotesleft Main Lemma'' of Moschovakis \cite{Mo74} (as in the proof of Theorem \ref{LemmacoINDGames}): since operators definable from $X$ by an open game quantifier are $X$-inductive (by Corollary \ref{TheoremCharacIndXGames}) and $\Game^\mathbb{R}\SIGMA^0_2$ is a Spector class, it suffices to show that $\Game^\mathbb{R}\SIGMA^0_2$ is closed under the inductive step. But the inductive step is positive analyical on $X$ and a real parameter, so clearly $\Game^\mathbb{R}\SIGMA^0_2$ is closed under it. 
\endproof
Write
\begin{equation*}
\phi'(s,X) = \text{\quotesleft Player I has a winning strategy in $H(X,s)$.''}
\end{equation*}
Clearly, $\phi'(s,X)$ is an operator in $\Game^\mathbb{R}\PI^0_1$ and thus in $\coind$ (by Corollary \ref{TheoremCharacIndXGames}). By Theorem \ref{theoremGeneralizedRichterAczel}, there is some $\eta\leq\sigma_\mathbb{R}$ such that $\phi'^\infty = \phi'^\eta$. By Corollary \ref{corollaryCharacterizationSigma02} and the hypothesis of the lemma,
$$L_{\sigma_\mathbb{R}}(\mathbb{R})\models\ad.$$
Since $\sigma_\mathbb{R}$ is recursively $\mathbb{R}$-inaccessible, i.e., $\mathbb{R}$-admissible and a limit of $\mathbb{R}$-admissibles, it follows that for every $\zeta<\eta$, every subset of $\mathbb{R}$ in $\Game^\mathbb{R}\SIGMA^0_1(\phi'^\zeta)$ is determined. Hence, an induction along $\zeta$ shows that
$$\phi'^\zeta = \phi^\zeta$$
for all $\zeta<\eta$. It follows that 
$$\phi^\infty = \phi'^\infty \in \Sigma_1^{L_{\sigma_\mathbb{R}}(\mathbb{R})} = \Game^\mathbb{R}\SIGMA^0_2.$$ 
(In case $\eta = \sigma^\mathbb{R}$, the last step of the induction follows from Sublemma \ref{SublemmaClosureInductive}.)

Now, suppose that $s \not\in \phi^\infty$. By definition, Player I does not have a winning strategy in $G(\phi^\infty,s)$. By the two sublemmata, $G(\phi^\infty,s)$ is determined. Thus, Player II has a winning strategy in $G(\phi^\infty,s)$. But then, an argument as in the proof of Theorem \ref{LemmacoINDGames} shows how to turn this into a winning strategy for $G^*(s)$. This concludes the proof of the lemma.
\endproof

We are now led to the main theorem of the article.

\begin{theorem}\label{TheoremMainFsigmaRestated}
The following are equivalent:
\begin{enumerate}
\item $\SIGMA^0_2$-determinacy for games of length $\omega^2$;
\item There is a $\Pi^+_1$-reflecting model of $\ad$ containing $\mathbb{R}$.
\end{enumerate}
\end{theorem}
\proof
If $\SIGMA^0_2$ games of length $\omega^2$ are determined, then clearly sets in $\Game^\mathbb{R}\SIGMA^0_2$ are determined, so that $L_{\sigma_\mathbb{R}}(\mathbb{R})\models\ad$, by Corollary \ref{corollaryCharacterizationSigma02}.

Conversely, suppose that there is a $\Pi^+_1$-reflecting model of $\ad$ containing $\mathbb{R}$. By Corollary \ref{CorollaryCharacPiplusRef}, $L_{\sigma_\mathbb{R}}(\mathbb{R})\models\ad$. Since $\sigma_\mathbb{R}$ is $\mathbb{R}$-admissible, 
\[\mathcal{P}(\mathbb{R})\cap L_{\sigma_\mathbb{R}}(\mathbb{R}) = \mathcal{P}(\mathbb{R})\cap \Delta_1^{L_{\sigma_\mathbb{R}}(\mathbb{R})}.\]
Let $\Gamma$ be the pointclass of all sets of reals which are $\Sigma_1$-definable over $L_{\sigma_\mathbb{R}}(\mathbb{R})$ with parameters.
Using that $L_{\sigma_\mathbb{R}}(\mathbb{R})\models\ad$, Steel \cite[Theorem 2.1]{St08} implies that $\Gamma$ has the scale property. Since $\sigma_\mathbb{R}$ is $\mathbb{R}$-admissible, $\Gamma$ is closed under existential and universal real quantification and thus by Moschovakis \cite{Mo71}, it has the uniformization property.
Hence, the hypotheses for the Kechris-Woodin determinacy transfer theorem \cite{KW83} are satisfied, and we may apply it to conclude that all sets in $\Sigma_1^{L_{\sigma_\mathbb{R}}(\mathbb{R})}$ are determined. By Corollary \ref{corollaryCharacterizationSigma02}, all sets in $\Game^\mathbb{R}\SIGMA^0_2$ are determined. By Lemma \ref{LemmaShorteningSigma02}, $\SIGMA^0_2$-games of length $\omega^2$ are determined.
\endproof

\section{Concluding Remarks}
Let us remark that the proof of Theorem \ref{TheoremMainFsigmaRestated} also yields the following result, which
was pointed out in \cite{Ag18a}: 
\begin{theorem}\label{TheoremKPOpen}
The following are equivalent over $\zfc$:
\begin{enumerate}
\item $\SIGMA^0_1$-determinacy for games of length $\omega^2$,
\item There is a transitive model of $\kp + \ad$ containing the reals.
\end{enumerate}
\end{theorem}
To prove Theorem \ref{TheoremKPOpen}, repeat the proof of Theorem \ref{TheoremMainFsigmaRestated}, using the least $\kappa$ such that $L_\kappa(\mathbb{R})$ is admissible in place of $\sigma_\mathbb{R}$, Corollary \ref{TheoremCharacIndXGames} in place of Corollary \ref{corollaryCharacterizationSigma02}, and Sublemma \ref{SublemmaClopen} in place of Lemma \ref{LemmaShorteningSigma02}.

A classical theorem of Kripke states that, in $L$, the least $\alpha$ such that $L_\alpha$ is a model of $\kp + \Sigma_1$-separation is the union of a chain of $\Sigma_1$-elementary substructures. The situation is analogous in $L(\mathbb{R})$.
This and Theorem \ref{mainintrofsigma} imply the consistency results mentioned in the introduction. 

\begin{theorem}\label{TheoremSigma02Consistency}
The following theories are strictly increasing in consistency strength:
\begin{enumerate}
\item \label{TheoremSigma02Consistency1} $\kp + \dc + \quotesleft \mathbb{R}$ exists'' $ + $ $\ad$;
\item \label{TheoremSigma02Consistency2} $\kp + \dc + \quotesleft \mathbb{R}$ exists'' $ + $ $\SIGMA^0_2$-determinacy for games of length $\omega^2$;
\item $\kp + \dc +\quotesleft \mathbb{R}$ exists'' $ + $ $\ad$ $ + $ $\Sigma_1$-separation.
\end{enumerate}
\end{theorem}
\proof
Working in $\kp + \dc + \quotesleft \mathbb{R}$ exists,'' assume that $\SIGMA^0_2$ games of length $\omega^2$ are determined. We show that there is a model of
\[\text{$\kp + \dc + \quotesleft \mathbb{R}$ exists'' $ + $ $\ad$.}\] 
Consider the following game of length $\omega^2$:
\begin{enumerate}
\item Player I begins by choosing a $\SIGMA^0_1$ game of length $\omega^2$;
\item Player II decides which role she wants to have in the game;
\item Players I and II play the game chosen by Player I, taking the roles specified by Player II.
\end{enumerate}
This game is perhaps not $\SIGMA^0_1$, but it is certainly, say, $\Delta^0_2$. It clearly cannot be won by Player I, and any winning strategy for Player II easily reduces continuously to a winning strategy for some player for any prescribed $\SIGMA^0_1$ game of length $\omega^2$. Since  $\mathbb{R}$ exists, one may use collection to conclude that there is a set $A$ containing winning strategies for all $\SIGMA^0_1$ games of length $\omega^2$. From $A$ and $\mathbb{R}$ one may use $\Sigma_0$-separation to conclude that there is a universal $\Game^\mathbb{R}\SIGMA^0_1$ set. Hence, there is a universal $\ind$ set, from which one can easily compute a prewellordering of $\mathbb{R}$ whose length is an $\mathbb{R}$-admissible ordinal $\kappa$. By collection, this ordinal exists, so $L_\kappa(\mathbb{R})\models\kp + \dc$ and
\[\ind = \mathcal{P}(\mathbb{R})\cap \Sigma^{L_\kappa(\mathbb{R})}_1.\]
Since all $\SIGMA^0_2$ games of length $\omega^2$ are determined, all $\Game^\mathbb{R}\SIGMA^0_1$ games of length $\omega$ are determined, and so 
\[L_\kappa(\mathbb{R})\models\ad.\]

For the remaining implication, suppose that $\ad$ holds and $\Sigma_1$-separation holds. It follows that 
\[L(\mathbb{R})\models \ad + \text{$\Sigma_1$-Separation}.\]
To see this, we assume without loss of generality that for no ordinal $\eta$ do we have 
\[L_\eta(\mathbb{R})\models \ad + \text{$\Sigma_1$-Separation},\]
i.e., that the least non-$\mathbb{R}$-projectible ordinal does not exist. Thus, it suffices to show that every subset of $\mathbb{R}$ which is $\Sigma_1$-definable (with parameters) over $L(\mathbb{R})$ belongs to $L(\mathbb{R})$. If not, then by Steel \cite[Lemma 1.14]{St08}, there is a partial surjection from $\mathbb{R}$ onto $L(\mathbb{R})$ which is $\Sigma_1$-definable over $L(\mathbb{R})$ with parameters. However, $L(\mathbb{R})$ is a $\Sigma_1$-definable class and $\Sigma_1$-separation holds (in $V$), so this surjection is actually a set. But this is impossible, for its restriction to the ordinals is also a surjection. Thus,
\[L(\mathbb{R})\models \ad + \text{$\Sigma_1$-Separation}\]
as claimed. Working in $L(\mathbb{R})$, there are arbitrarily large $\mathbb{R}$-stable ordinals. But it is immediate from the definition that if
\[L_\alpha(\mathbb{R})\prec_1 L_{\alpha^++1}(\mathbb{R}),\]
(recall Definition \ref{DefinitionPlusLR}) then $\sigma_\mathbb{R}<\alpha$, so there is a (set) model of $\kp$ in which all games of length $\omega^2$ are determined (and $\ad$ holds in this model, additionally).
\endproof

By a slightly more elaborate argument involving Inner Model Theory (one can e.g., use the results in Part II of M\"uller \cite{Uh16}), one can replace $\kp$ by $\zfc$ in the statement of Theorem \ref{TheoremSigma02Consistency}\eqref{TheoremSigma02Consistency2} (and, in fact, by much stronger theories). Similarly, and using Theorem \ref{TheoremKPOpen}, one can replace Theorem \ref{TheoremSigma02Consistency}\eqref{TheoremSigma02Consistency1} by the theory $\zfc$ + $\SIGMA^0_1$-determinacy for games of length $\omega^2$. We omit the details.

\bibliographystyle{abbrv}
\bibliography{References}

\end{document}